# Pile-up probabilities for the Laplace likelihood estimator of a non-invertible first order moving average

F. Jay Breidt[1,*,†], Richard A. Davis[1,†,‡], Nan-Jung Hsu[2] and Murray Rosenblatt[3]

*Colorado State University, National Tsing-Hua University and University of California at San Diego*

**Abstract:** The first-order moving average model or MA(1) is given by $X_t = Z_t - \theta_0 Z_{t-1}$, with independent and identically distributed $\{Z_t\}$. This is arguably the simplest time series model that one can write down. The MA(1) with unit root ($\theta_0 = 1$) arises naturally in a variety of time series applications. For example, if an underlying time series consists of a linear trend plus white noise errors, then the differenced series is an MA(1) with unit root. In such cases, testing for a unit root of the differenced series is equivalent to testing the adequacy of the trend plus noise model. The unit root problem also arises naturally in a signal plus noise model in which the signal is modeled as a random walk. The differenced series follows a MA(1) model and has a unit root if and only if the random walk signal is in fact a constant.

The asymptotic theory of various estimators based on Gaussian likelihood has been developed for the unit root case and nearly unit root case ($\theta = 1 + \beta/n, \beta \leq 0$). Unlike standard $1/\sqrt{n}$-asymptotics, these estimation procedures have $1/n$-asymptotics and a so-called pile-up effect, in which $P(\hat{\theta} = 1)$ converges to a positive value. One explanation for this pile-up phenomenon is the lack of identifiability of $\theta$ in the Gaussian case. That is, the Gaussian likelihood has the same value for the two sets of parameter values $(\theta, \sigma^2)$ and $(1/\theta, \theta^2\sigma^2)$. It follows that $\theta = 1$ is always a critical point of the likelihood function. In contrast, for non-Gaussian noise, $\theta$ is identifiable for all real values. Hence it is no longer clear whether or not the same pile-up phenomenon will persist in the non-Gaussian case. In this paper, we focus on limiting pile-up probabilities for estimates of $\theta_0$ based on a Laplace likelihood. In some cases, these estimates can be viewed as Least Absolute Deviation (LAD) estimates. Simulation results illustrate the limit theory.

## 1. Introduction

The moving average model of order one (MA(1)) given by

$$(1.1) \qquad X_t = Z_t - \theta_0 Z_{t-1},$$

---

[1]Department of Statistics, Colorado State University, Ft. Collins, CO 80523, USA, e-mail: `jbreidt@stat.colostate.edu`; `rdavis@stat.colostate.edu`
[2]Institute of Statistics, National Tsing-Hua University, Hsinchu, Taiwan, e-mail: `njhsu@stat.nthu.edu.tw`
[3]Department of Mathematics, University of California, San Diego, La Jolla, CA 92093, USA, e-mail: `mrosenblatt@ucsd.edu`
*Research supported by NSF grant DMS-9972015.
†Research supported by EPA STAR grant CR-829095.
‡Research supported by NSF grant DMS-0308109.

where $\{Z_t\}$ is a sequence of independent and identically distributed random variables with mean 0 and variance $\sigma^2$, is one of the simplest models in time series. The MA(1) model is invertible if and only if $|\theta_0| < 1$, since in this case $Z_t$ can be represented explicitly in terms of past values of the $X_t$, i.e.,

$$Z_t = \sum_{j=0}^{\infty} \theta_0^j X_{t-j}.$$

Under this invertibility constraint, standard estimation procedures that produce asymptotically normal estimates are readily available. For example, if $\hat{\theta}$ represents the maximum likelihood estimator, found by maximizing the Gaussian likelihood based on the data $X_1, \ldots, X_n$, then it is well known (see Brockwell and Davis [3]), that

(1.2) $$\sqrt{n}(\hat{\theta} - \theta_0) \xrightarrow{d} N(0, 1 - \theta_0^2).$$

From the form of the limiting variance in (1.2), the asymptotic behavior of $\hat{\theta}$, let alone the scaling, is not immediately clear in the unit root case corresponding to $\theta_0 = 1$.

In the Gaussian case, the parameters $\theta_0$ and $\sigma^2$ are not identifiable without the constraint $|\theta_0| \leq 1$. In particular, the profile Gaussian log-likelihood, obtained by concentrating out the variance parameter, satisfies

$$L(\theta) = L(1/\theta).$$

It follows that $\theta = 1$ is a critical value of the profile likelihood and hence there is a positive probability that $\theta = 1$ is indeed the maximum likelihood estimator. If $\theta_0 = 1$, then it turns out that this probability does not vanish asymptotically (see for example Anderson and Takemura [1], Tanaka [7], and Davis and Dunsmuir [6]). This phenomenon is referred to as the pile-up effect. For the case that $\theta_0 = 1$ or is near one in the sense that $\theta_0 = 1 + \gamma/n$, it was shown in Davis and Dunsmuir [6] that

$$n(\hat{\theta} - \theta_0) \xrightarrow{d} \xi_\gamma,$$

where $\xi_\gamma$ is random variable with a discrete component at 0, corresponding to the asymptotic pile-up probability, and a continuous component on $(-\infty, 0)$.

The MA(1) with unit root ($\theta_0 = 1$) arises naturally in a variety of time series applications. For example, if an underlying time series consists of a linear trend plus white noise errors, then the differenced series is an MA(1) with a unit root. In such cases, testing for a unit root of the differenced series is equivalent to testing the adequacy of the trend plus noise model. The unit root problem also arises naturally in a signal plus noise model in which the signal is modeled as a random walk. The differenced series follows a MA(1) model and has a unit root if and only if the random walk signal is in fact a constant.

For Gaussian likelihood estimation, the pile-up effect is directly attributable to the non-identifiability of $\theta_0$ in the unconstrained parameter space. On the other hand, if the data are non-Gaussian, then $\theta_0$ is identifiable (see Breidt and Davis [2]). In this paper, we focus on the pile-up probability for estimates based on a Laplace likelihood. Assuming a Laplace distribution for the noise, we derive an expression for the joint likelihood of $\theta$ and $z_{init}$, where $z_{init}$ is an augmented variable that is treated as a parameter and the scale parameter $\sigma$ is concentrated out of the likelihood. If $z_{init}$ is set equal to 0, then the resulting joint likelihood corresponds to the least absolute deviation (LAD) objective function and the estimator of $\theta$ is referred to as the LAD estimator of $\theta_0$. The exact likelihood can be obtained by integrating out $z_{init}$. In this case the resulting estimator is referred to as the quasi-maximum likelihood estimator of $\theta_0$. It turns out that the estimator based on maximizing the joint likelihood always has a positive pile-up probability in the limit regardless of the true noise distribution. In contrast, the quasi-maximum likelihood estimator has a limiting pile-up probability of zero.

In Section 2, we describe the main asymptotic results. We begin by deriving an expression for computing the joint likelihood function based on the observed data and the augmented variable $Z_{init}$, in terms of the density function of the noise. The exact likelihood function can then be computed by integrating out $Z_{init}$. After a reparameterizion, we derive the limiting behavior of the joint likelihood for the case when the noise is assumed to follow a Laplace distribution. In Section 3, we focus on the problem of calculating asymptotic pile-up probabilities for estimators which minimize the joint Laplace likelihood (as a function of $\theta$ and $z_{init}$) and the exact Laplace likelihood. Section 4 contains simulation results which illustrate the asymptotic theory of Section 3.

## 2. Main result

Let $\{X_t\}$ be the MA(1) model given in (1.1) where $\theta_0 \in \mathbb{R}$, $\{Z_t\}$ is a sequence of iid random variables with $EZ_t = 0$ and density function $f_Z$. In order to compute the likelihood based on the observed data $\boldsymbol{X}_n = (X_1, \ldots, X_n)'$, it is convenient to define an augmented initial variable $Z_{init}$ defined by

$$Z_{init} = \begin{cases} Z_0, & \text{if } |\theta| \leq 1, \\ Z_n - \sum_{t=1}^n X_t, & \text{otherwise.} \end{cases}$$

A straightforward calculation shows that the joint density of the observed data $\boldsymbol{X}_n = (X_1, X_2, \ldots, X_n)'$ and the initial variable $Z_{init}$ satisfies

$$f_{\boldsymbol{X}, Z_{init}}(\boldsymbol{x}_n, z_{init}) = \prod_{j=0}^n f_Z(z_j) \left(1_{\{|\theta| \leq 1\}} + |\theta|^{-n} 1_{\{|\theta| > 1\}}\right),$$

where the *residuals* $\{z_t\}$ are functions of $\boldsymbol{X}_n = \boldsymbol{x}_n$, $\theta$, and $Z_{init} = z_{init}$ which can be solved forward by $z_t = X_t + \theta z_{t-1}$ for $t = 1, 2, \ldots, n$ with the initial $z_0 = z_{init}$ if $|\theta| \leq 1$ and backward by $z_{t-1} = \theta^{-1}(z_t - X_t)$ for $t = n, n-1, \ldots, 1$ with the initial $z_n = z_{init} + \sum_{t=1}^n X_t$, if $|\theta| > 1$.

The Laplace log-likelihood is obtained by taking the density function for $Z_t$ to be $f_Z(z) = \exp\{-|z|/\sigma\}/(2\sigma)$. If we view $z_{init}$ as a parameter, then the *joint* log-likelihood is given by

$$(2.1) \qquad -(n+1)\log 2\sigma - \frac{1}{\sigma} \sum_{t=0}^n |z_t| - n(\log |\theta|) 1_{\{|\theta| > 1\}}.$$

Maximizing this function with respect to the scale parameter $\sigma$, we obtain

$$\hat{\sigma} = \sum_{t=0}^n |z_t|/(n+1).$$

It follows that maximizing the joint Laplace log-likelihood is equivalent to minimizing the following objective function,

$$\ell_n(\theta, z_{init}) = \begin{cases} \sum_{t=0}^{n} |z_t|, & \text{if } |\theta| \leq 1, \\ \sum_{t=0}^{n} |z_t||\theta|, & \text{otherwise.} \end{cases} \tag{2.2}$$

In order to study the asymptotic properties of the minimizer of $\ell_n$ when the model $\theta_0 = 1$, we follow Davis and Dunsmuir [6] by building the sample size into the parameterization of $\theta$. Specifically, we use

$$\theta = 1 + \frac{\beta}{n}, \tag{2.3}$$

where $\beta$ is any real number. Additionally, since we are also treating $z_{init}$ as a parameter, this term is reparameterized as

$$z_{init} = Z_0 + \frac{\alpha\sigma}{\sqrt{n}}. \tag{2.4}$$

Under the $(\beta, \alpha)$ parameterization, minimizing $\ell_n$ with respect to $\theta$ and $z_{init}$ is equivalent to minimizing the function,

$$U_n(\beta, \alpha) \equiv \frac{1}{\sigma}\left[\ell_n(\theta, z_{init}) - \ell_n(1, Z_0)\right],$$

with respect to $\beta$ and $\alpha$. The following theorem describes the limiting behavior of $U_n$.

**Theorem 2.1.** *For the model (1.1) with $\theta_0 = 1$, assume the noise sequence $\{Z_t\}$ is IID with $EZ_t = 0$, $E[\text{ sign}(Z_t)] = 0$ (i.e., median of $Z_t$ is zero), $EZ_t^4 < \infty$ and common probability density function $f_Z(z) = \sigma^{-1}f(z/\sigma)$, where $\sigma > 0$ is the scale parameter. We further assume that the density function $f_Z$ has been normalized so that $\sigma = E|Z_t|$. Then*

$$U_n(\beta, \alpha) \xrightarrow{fidi} U(\beta, \alpha), \tag{2.5}$$

*where $\xrightarrow{fidi}$ denotes convergence in distribution of finite dimensional distributions and*

$$U(\beta, \alpha) = \int_0^1 \left[\beta \int_0^s e^{\beta(s-t)} dS(t) + \alpha e^{\beta s}\right] dW(s)$$
$$+ f(0) \int_0^1 \left[\beta \int_0^s e^{\beta(s-t)} dS(t) + \alpha e^{\beta s}\right]^2 ds, \tag{2.6}$$

*for $\beta \leq 0$, and*

$$U(\beta, \alpha) = \int_0^1 \left[-\beta \int_{s+}^1 e^{-\beta(t-s)} dS(t) + \alpha e^{-\beta(1-s)}\right] dW(s)$$
$$+ f(0) \int_0^1 \left[-\beta \int_s^1 e^{-\beta(t-s)} dS(t) + \alpha e^{-\beta(1-s)}\right]^2 ds, \tag{2.7}$$

*for $\beta > 0$, in which $S(t)$ and $W(t)$ are the limits of the following partial sums*

$$S_n(t) = \frac{1}{\sqrt{n}} \sum_{i=0}^{[nt]} Z_i/\sigma, \quad W_n(t) = \frac{1}{\sqrt{n}} \sum_{i=0}^{[nt]} \text{sign}(Z_i),$$

*respectively.*

**Remark.** The stochastic integrals in (2.6) and (2.7) refer to Itô integrals. The double stochastic stochastic integral in the first term on the right side of (2.7) is computed as

$$\int_0^1 \int_{s+}^1 e^{-\beta(t-s)} dS(t) dW(s) = \int_0^1 e^{-\beta t} dS(t) \int_0^1 e^{\beta s} dW(s)$$
$$- \int_0^1 \int_0^s e^{-\beta(t-s)} dS(t) dW(s) - \int_0^1 dS(t) dW(t),$$

where (see (2.15) below)

$$\int_0^1 dS(t) dW(t) = E(Z_i \, sign(Z_i))/\sigma = E|Z_i|/\sigma = 1 \, .$$

*Proof.* We only prove the result (2.5) for a fixed $(\beta, \alpha)$; the extension to a finite collection of $(\beta, \alpha)$'s is relatively straightforward. First consider the case $\beta \leq 0$. For calculating the Laplace likelihood $\ell_n(\theta, z_{init})$ based on model (1.1), the residuals are solved by $z_t = X_t + \theta z_{t-1}$ for $t = 1, 2, \ldots, n$ with the initial value $z_0 = z_{init}$. Since $X_t = Z_t - Z_{t-1}$, all of the true innovations can be solved forward by $Z_t = X_t + Z_{t-1}$ for $t = 1, 2, \ldots, n$ with the initial $Z_0$. Therefore, the centered term $\ell_n(1, Z_0)$ can be written as

$$\ell_n(1, Z_0) = |Z_0| + \sum_{i=1}^n |X_i + X_{i-1} + \cdots + X_1 + Z_0| = \sum_{i=0}^n |Z_i|.$$

For $\beta \leq 0$, i.e., $\theta \leq 1$,

$$z_i = X_i + \theta X_{i-1} + \cdots + \theta^{i-1} X_1 + \theta^i z_{init}$$
$$= (Z_i - Z_{i-1}) + \theta(Z_{i-1} - Z_{i-2}) + \cdots + \theta^{i-1}(Z_1 - Z_0) + \theta^i z_{init}$$
$$= Z_i - (1-\theta) Z_{i-1} - \theta(1-\theta) Z_{i-2} - \cdots - \theta^{i-1}(1-\theta) Z_0 - \theta^i (Z_0 - z_{init}),$$

which, under the true model $\theta = 1$, implies

(2.8)
$$\frac{1}{\sigma} [\ell_n(\theta, z_{init}) - \ell_n(1, Z_0)] = \frac{1}{\sigma} \left( \sum_{i=0}^n |z_i| - \sum_{i=0}^n |Z_i| \right)$$
$$= \frac{1}{\sigma} \sum_{i=0}^n (|Z_i - y_i| - |Z_i|),$$

where $y_0 \equiv Z_0 - z_{init}$ and

$$y_i \equiv (1-\theta) \sum_{j=0}^{i-1} \theta^{i-1-j} Z_j + \theta^i (Z_0 - z_{init}),$$

for $i = 1, 2, \ldots, n$. Using the identity

(2.9) $\quad |Z - y| - |Z| = -y \, sign(Z) + 2(y - Z) \left( 1_{\{0 < Z < y\}} - 1_{\{y < Z < 0\}} \right)$

for $Z \neq 0$, the equation (2.8) is expressed as two summations, the first of which is

$$
\begin{aligned}
-\sum_{i=0}^{n} \frac{y_i}{\sigma} \operatorname{sign}(Z_i) &= (\theta - 1) \sum_{i=1}^{n} \left( \sum_{j=0}^{i-1} \theta^{i-1-j} \frac{Z_j}{\sigma} \right) \operatorname{sign}(Z_i) \\
&\quad + \frac{z_{init} - Z_0}{\sigma} \sum_{i=0}^{n} \theta^i \operatorname{sign}(Z_i) \\
&= \frac{\beta}{n} \sum_{i=1}^{n} \left[ \sum_{j=0}^{i-1} \left(1 + \frac{\beta}{n}\right)^{i-j-1} \frac{Z_j}{\sigma} \right] \operatorname{sign}(Z_i) \\
&\quad + \frac{\alpha}{\sqrt{n}} \sum_{i=0}^{n} \left(1 + \frac{\beta}{n}\right)^i \operatorname{sign}(Z_i) \\
&= \beta \int_0^1 \int_0^{s-} \left(1 + \frac{\beta}{n}\right)^{-nt} dS_n(t) \left(1 + \frac{\beta}{n}\right)^{ns-1} dW_n(s) \\
&\quad + \alpha \int_0^1 \left(1 + \frac{\beta}{n}\right)^{ns} dW_n(s) \\
&\to \beta \int_0^1 \int_0^s e^{\beta(s-t)} dS(t) dW(s) + \alpha \int_0^1 e^{\beta s} dW(s),
\end{aligned}
$$
(2.10)

where the limit in (2.10) follows from a simple adaptation of Theorem 2.4 (ii) in Chan and Wei [4].

To handle the second summation in computing $U_n(\beta, \alpha)$, we approximate the sum

$$\sum_{i=0}^{n} 2 \frac{y_i - Z_i}{\sigma} \left( 1_{\{0 < Z_i < y_i\}} - 1_{\{y_i < Z_i < 0\}} \right)$$

by

$$\sum_{i=0}^{n} 2E \left[ \frac{y_i - Z_i}{\sigma} \left( 1_{\{0 < Z_i < y_i\}} - 1_{\{y_i < Z_i < 0\}} \right) | \mathcal{F}_{i-1} \right],$$

where $\mathcal{F}_i$ is the $\sigma$-field generated by $\{Z_j : j = 0, 1, \ldots, i\}$. First we establish convergence of the latter sum and then show that the variance of the difference in sums converges to zero. Since

$$\max_{1 \leq i \leq n} |y_i| \to 0,$$

$y_i \in \mathcal{F}_{i-1}$, we have

$$
\begin{aligned}
2E \left[ \left( \frac{y_i - Z_i}{\sigma} \right) 1_{\{0 < Z_i < y_i\}} | \mathcal{F}_{i-1} \right] &= 2 \int_0^{y_i} \left( \frac{y_i - Z}{\sigma} \right) \frac{1}{\sigma} f(\frac{z}{\sigma}) dz \\
&\approx f(0) \int_0^{y_i} 2 \left( \frac{y_i - z}{\sigma} \right) d \left( \frac{z}{\sigma} \right) \\
&= f(0) \left( \frac{y_i}{\sigma} \right)^2,
\end{aligned}
$$

for $y_i > 0$, and

$$2E\left[\left(\frac{y_i - Z_i}{\sigma}\right) 1_{\{y_i < Z_i < 0\}} | \mathcal{F}_{i-1}\right] = 2\int_{y_i}^0 \left(\frac{y_i - z}{\sigma}\right) \frac{1}{\sigma} f(\frac{z}{\sigma}) dz$$
$$\approx f(0) \int_{y_i}^0 2\left(\frac{y_i - z}{\sigma}\right) d\left(\frac{z}{\sigma}\right)$$
$$= -f(0) \left(\frac{y_i}{\sigma}\right)^2,$$

for $y_i < 0$. Combining these two cases, we have

$$2\sum_{i=0}^n E\left[\frac{y_i - Z_i}{\sigma} \left(1_{\{0 < Z_i < y_i\}} - 1_{\{y_i < Z_i < 0\}}\right) | \mathcal{F}_{i-1}\right] \approx f(0) \sum_{i=0}^n \left(\frac{y_i}{\sigma}\right)^2,$$

where

(2.11)
$$\sum_{i=0}^n \left(\frac{y_i}{\sigma}\right)^2 = \sum_{i=0}^n \left\{(1-\theta) \sum_{j=1}^{i-1} \theta^{i-1-j} \frac{Z_j}{\sigma} + \theta^i \frac{Z_0 - z_0}{\sigma}\right\}^2$$
$$= \sum_{i=1}^n \left[\frac{-\beta}{n} \sum_{j=1}^{i-1} \left(1 + \frac{\beta}{n}\right)^{i-1-j} \frac{Z_j}{\sigma} - \frac{\alpha}{\sqrt{n}} \left(1 + \frac{\beta}{n}\right)^i\right]^2$$
$$= \sum_{i=1}^n \left[\beta \int_0^{(i-1)/n} \left(1 + \frac{\beta}{n}\right)^{i-1-sn} dS_n(s) + \alpha \left(1 + \frac{\beta}{n}\right)^i\right]^2 \frac{1}{n}$$
$$\to \int_0^1 \left[\beta \int_0^s e^{\beta(s-t)} dS(t) + \alpha e^{\beta s}\right]^2 ds$$

in distribution as $n \to \infty$.

It is left to show that

(2.12)
$$2\sum_{i=0}^n \frac{y_i - Z_i}{\sigma} \left(1_{\{0 < Z_i < y_i\}} - 1_{\{y_i < Z_i < 0\}}\right)$$
$$- 2\sum_{i=0}^n E\left[\frac{y_i - Z_i}{\sigma} \left(1_{\{0 < Z_i < y_i\}} - 1_{\{y_i < Z_i < 0\}}\right) | \mathcal{F}_{i-1}\right]$$

converges to zero in probability. Define

$$y_i^* \equiv 2\frac{y_i - Z_i}{\sigma} \left(1_{\{0 < Z_i < y_i\}} - 1_{\{y_i < Z_i < 0\}}\right).$$

The expectation of (2.12) is zero and therefore, it is enough to show that the variance of (2.12) also converges to zero. The variance of (2.12) is equal to

$$\sum_{i=0}^{n} \operatorname{var}(y_i^* - E(y_i^*|\mathcal{F}_{i-1})) + 2 \sum_{i<j} \operatorname{cov}\left(y_i^* - E(y_i^*|\mathcal{F}_{i-1}), y_j^* - E(y_j^*|\mathcal{F}_{j-1})\right)$$

$$= \sum_{i=0}^{n} E\left[y_i^* - E(y_i^*|\mathcal{F}_{i-1})\right]^2$$

$$= \sum_{i=0}^{n} EE\left[(y_i^*)^2 - (E(y_i^*|\mathcal{F}_{i-1}))^2 \mid \mathcal{F}_{i-1}\right]$$

(2.13)
$$= \sum_{i=0}^{n} E\left[E\left((y_i^*)^2|\mathcal{F}_{i-1}\right) - (E(y_i^*|\mathcal{F}_{i-1}))^2\right]$$

$$\approx \sum_{i=0}^{n} E\left[\frac{4}{3}f(0)\left(\frac{y_i}{\sigma}\right)^3 - f(0)^2 \left(\frac{y_i}{\sigma}\right)^4\right]$$

$$\approx \frac{4}{3}f(0)E\left[\sum_{i=0}^{n}\left(\frac{y_i}{\sigma}\right)^3\right] - f(0)^2 E\left[\sum_{i=0}^{n}\left(\frac{y_i}{\sigma}\right)^4\right]$$

$$\to 0,$$

as $n \to \infty$, where

$$\operatorname{cov}\left(y_i^* - E(y_i^*|\mathcal{F}_{i-1}), y_j^* - E(y_j^*|\mathcal{F}_{j-1})\right)$$
$$= E\left[y_i^* - E(y_i^*|\mathcal{F}_{i-1})\right]\left[y_j^* - E(y_j^*|\mathcal{F}_{j-1})\right]$$
$$= EE\left[(y_i^* - E(y_i^*|\mathcal{F}_{i-1}))\left(y_j^* - E(y_j^*|\mathcal{F}_{j-1})\right) \Big| \mathcal{F}_{j-1}\right]$$
$$= E\left[(y_i^* - E(y_i^*|\mathcal{F}_{i-1}))E\left(y_j^* - E(y_j^*|\mathcal{F}_{j-1}) \Big| \mathcal{F}_{j-1}\right)\right]$$
$$= 0,$$

for $i < j$, and

$$E(y_i^*|\mathcal{F}_{i-1}) \approx f(0)\left(\frac{y_i}{\sigma}\right)^2,$$

$$E\left((y_i^*)^2|\mathcal{F}_{i-1}\right) \approx \frac{4}{3}f(0)\left(\frac{y_i}{\sigma}\right)^3,$$

$$\sqrt{n}\sum_{i=0}^{n}\left(\frac{y_i}{\sigma}\right)^3 \to -\int_0^1 \left(\beta \int_0^s e^{\beta(s-t)}dS(t) + \alpha e^{\beta s}\right)^3 ds,$$

$$n\sum_{i=0}^{n}\left(\frac{y_i}{\sigma}\right)^4 \to \int_0^1 \left(\beta \int_0^s e^{\beta(s-t)}dS(t) + \alpha e^{\beta s}\right)^4 ds.$$

Based on (2.10), (2.11), and (2.13), the proof for $\beta \leq 0$ is complete.

The proof for $\beta \geq 0$ given in (2.7) is similar to that for $\beta \leq 0$. For $\beta \geq 0$, i.e., $\theta \geq 1$, the residuals $\{z_t\}$ are solved backward by $z_{t-1} = \theta^{-1}(z_t - X_t)$ for $t = n, n-1, \ldots, 1$ with the initial $z_n \equiv z_{init} + \sum_{t=1}^{n} X_t$. Solving these equations, we have

$$z_{n-1-i} = -\theta^{-1}\left(X_{n-i} + \theta^{-1}X_{n-i-1} + \cdots + \theta^{-i}X_n - \theta^{-i}z_n\right),$$

for $i = 0, 1, \ldots, n-1$. Writing $X_t = Z_t - Z_{t-1}$, we obtain

$$\begin{aligned}
-z_{n-1-i}\theta &= X_{n-i} + \theta^{-1}X_{n-i-1} + \cdots + \theta^{-i}X_n - \theta^{-i}z_n \\
&= (Z_{n-i} - Z_{n-i-1}) + \theta^{-1}(Z_{n-i+1} - Z_{n-i}) + \cdots \\
&\quad + \theta^{-i}(Z_n - Z_{n-1}) - \theta^{-i}z_n \\
&= -Z_{n-i-1} + (1 - \theta^{-1})Z_{n-i} + \cdots + \theta^{-(i-1)}(1 - \theta^{-1})Z_{n-1} \\
&\quad + \theta^{-i}(Z_n - z_n) \\
&= -Z_{n-i-1} + y_{n-i-1},
\end{aligned}$$

where

$$\begin{aligned}
y_{n-1-i} &\equiv (1 - \theta^{-1}) \sum_{j=1}^{i} (\theta^{-1})^{i-j} Z_{n-j} + \theta^{-i}(Z_n - z_n) \\
&= (1 - \theta^{-1}) \sum_{j=1}^{i} (\theta^{-1})^{i-j} Z_{n-j} + \theta^{-i}\left[\left(\sum_{i=1}^{n} X_i + Z_0\right) - \left(\sum_{i=1}^{n} X_i + z_{init}\right)\right] \\
&= (1 - \theta^{-1}) \sum_{j=1}^{i} (\theta^{-1})^{i-j} Z_{n-j} + \theta^{-i}(Z_0 - z_{init}),
\end{aligned}$$

for $i = 0, 1, \ldots, n-1$ and $y_n \equiv Z_n - z_n = Z_0 - z_{init}$. Again, for $\theta \geq 1$, we have

$$\frac{1}{\sigma}[\ell_n(\theta, z_{init}) - \ell_n(1, Z_0)] = \frac{1}{\sigma} \sum_{i=0}^{n} (|Z_i - y_i| - |Z_i|),$$

which has the same form as that for $\theta \leq 1$ but with different $\{y_i\}$. Following a similar derivation for $\theta \leq 1$, one can show that

$$-\sum_{i=1}^{n} \frac{y_i}{\sigma} \operatorname{sign}(Z_i) \to -\beta \int_0^1 \int_{s+}^1 e^{-\beta(t-s)} dS(t) dW(s) + \alpha \int_0^1 e^{-\beta(1-s)} dW(s),$$

$$\sum_{i=0}^{n} \frac{y_i^2}{\sigma^2} \to \int_0^1 \left[-\beta \int_s^1 e^{-\beta(t-s)} dS(t) + \alpha e^{-\beta(1-s)}\right]^2 ds,$$

in distribution as $n \to \infty$. Combining this with the analogous result (2.13) for $\beta \geq 0$, completes the proof. □

We close this section with some elementary results concerning the relationship between the limiting Brownian motions $S(t)$ and $W(t)$ that will be used in the sequel. Since $\sigma = E|Z_t|$, the process $S(t)$ can be decomposed as

(2.14) $$S(t) = W(t) + cV(t),$$

where $\{W(t)\}$ and $\{V(t)\}$ are independent standard Bronwnian motions on $[0, 1]$ and

$$c = \sqrt{\operatorname{Var}(Z_t)/\sigma^2 - 1}.$$

In addition, we have the following identities

$$\int_0^1 V(s)ds = V(1) - \int_0^1 sdV(s),$$
$$\int_0^1 V(s)dW(s) = V(1)W(1) - \int_0^1 W(s)dV(s),$$
$$\int_0^1 dW(s)dW(s) = \int_0^1 ds = 1,$$
$$\int_0^1 dV(s)dW(s) = 0,$$

where the first two equations can be obtained easily by integration by parts. It follows that

$$(2.15) \qquad \int_0^1 dS(s)dW(s) = \int_0^1 dW(s)dW(s) + c\int_0^1 dV(s)dW(s) = 1\,.$$

## 3. Pile-up probabilities

### 3.1. Joint likelihood

In this section, we will consider the local maximizer of the joint likelihood given by $-\ell_n$ in (2.2). This estimator was also studied by Davis and Dunsmuir [6] in the Gaussian case. Denote by $(\hat{\theta}_n^{(J)}, \hat{z}_{init,n}^{(J)})$ the local minimizer of $\ell_n(\theta, z_{init})$ in which $\hat{\theta}_n^{(J)}$ is closest to 1. Using the $(\beta, \alpha)$ parameterization given in (2.3) and (2.4), this is equivalent to finding the local minimizer $(\hat{\beta}_n^{(J)}, \hat{\alpha}_n^{(J)})$ of $U_n(\beta, \alpha)$ in which $\hat{\beta}_n^{(J)}$ is closest to zero. Moreover, the respective local minimizers of $\ell_n$ and $U_n$ are connected through the following relations:

$$(3.1) \qquad \hat{\theta}_n^{(J)} = 1 + \frac{\hat{\beta}_n^{(J)}}{n}, \quad \hat{z}_{init,n}^{(J)} = Z_0 + \frac{\hat{\alpha}_n^{(J)}\sigma}{\sqrt{n}}.$$

If the convergence of $U_n$ to $U$ in Theorem 1 is strengthened to weak convergence of processes on $C(\mathbb{R}^2)$, then the argument given in Davis and Dunsmuir [6] suggests the convergence in distribution of $(\hat{\beta}_n^{(J)}, \hat{\alpha}_n^{(J)})$ to $(\beta^{(J)}, \alpha^{(J)})$, where $(\hat{\beta}^{(J)}, \hat{\alpha}^{(J)})$ is the local minimizer of $U(\beta, \alpha)$ in which $\hat{\beta}^{(J)}$ is closest to 0. It follows that

$$(3.2) \qquad (n(\hat{\theta}_n^{(J)} - 1), \sqrt{n}(\hat{z}_{init,n}^{(J)} - Z_0)/\sigma) \xrightarrow{d} (\hat{\beta}^{(J)}, \hat{\alpha}^{(J)})\,.$$

The proofs of these results are the subject of on-going research and will appear in a forthcoming manuscript.

Turning to the question of pile-up probabilities, we have that 1 is a local minimizer if the derivative of the criterion function from the left is negative and the derivative from the right is positive; that is,

$$P(\hat{\theta}_n^{(J)} = 1) = P(\hat{\beta}_n^{(J)} = 0)$$
$$= P\left[\lim_{\beta\uparrow 0}\frac{\partial}{\partial\beta}U_n\left(\beta, \hat{\alpha}_n(\beta)\right) < 0 \text{ and } \lim_{\beta\downarrow 0}\frac{\partial}{\partial\beta}U_n\left(\beta, \hat{\alpha}_n(\beta)\right) > 0\right],$$

where $\hat{\alpha}_n(\beta) = \arg\min_\alpha U_n(\beta, \alpha)$ for given $\beta$. Assuming convergence of the right- and left-hand derivatives of the process $U_n(\beta, \hat{\alpha}_n(\beta))$, we obtain
(3.3)
$$\lim_{n\to\infty} P(\hat{\theta}_n^{(J)} = 1) = P\left[\lim_{\beta\uparrow 0} \frac{\partial}{\partial\beta} U(\beta, \hat{\alpha}(\beta)) < 0 \text{ and } \lim_{\beta\downarrow 0} \frac{\partial}{\partial\beta} U(\beta, \hat{\alpha}(\beta)) > 0\right],$$

where $\hat{\alpha}(\beta) = \arg\min_\alpha U(\beta, \alpha)$. We now proceed to simplify the limits of the two derivatives in the brackets of (3.3) in terms of the processes $S(t)$ and $W(t)$. According to (2.6) in Theorem 2.1, we have

$$\lim_{\beta\uparrow 0} \frac{\partial}{\partial\alpha} U(\beta, \alpha) = \lim_{\beta\uparrow 0} \left\{ \int_0^1 e^{\beta s} dW(s) + f(0) 2\alpha \int_0^1 e^{2\beta s} ds \right\}$$
$$= \int_0^1 dW(s) + 2\alpha f(0) \int_0^1 ds$$
$$= W(1) + 2\alpha f(0),$$

and therefore

$$\hat{\alpha}(0-) = -\frac{W(1)}{2f(0)}.$$

The derivative of $U(\beta, \alpha)$ with respect to $\beta$ at zero from the left-hand side satisfies

$$\frac{\partial}{\partial\beta} U(\beta, \alpha) = \int_0^1 \int_0^s e^{\beta(s-t)} dS(t) dW(s) + \beta \int_0^1 \int_0^s e^{\beta(s-t)}(s-t) dS(t) dW(s)$$
$$+ \alpha \int_0^1 e^{\beta s} s\, dW(s)$$
$$+ f(0) \left\{ 2\beta \int_0^1 \left( \int_0^s e^{\beta(s-t)} dS(t) \right)^2 ds \right.$$
$$+ \beta^2 \int_0^1 2 \left( \int_0^s e^{\beta(s-t)} dS(t) \right) \left( \int_0^s e^{\beta(s-t)}(s-t) dS(t) \right) ds$$
$$+ \alpha^2 \int_0^1 e^{2\beta s} 2s\, ds + 2\alpha \int_0^1 e^{\beta s} \left( \int_0^s e^{\beta(s-t)} dS(t) \right) ds$$
$$\left. + 2\alpha\beta \int_0^1 e^{\beta s} \left( \int_0^s e^{\beta(s-t)}(2s-t) dS(t) \right) ds \right\}.$$

Taking the limit as $\beta \uparrow 0$, we have

$$\lim_{\beta\uparrow 0} \frac{\partial}{\partial\beta} U(\beta, \hat{\alpha}(\beta)) = \int_0^1 \int_0^s dS(t) dW(s) + \hat{\alpha}(0-) \int_0^1 s\, dW(s)$$
$$+ f(0) \left\{ \hat{\alpha}^2(0-) \int_0^1 2s\, ds + 2\hat{\alpha}(0-) \int_0^1 \int_0^s dS(t) ds \right\}$$
(3.4)
$$= \int_0^1 S(s) dW(s) - W(1) \int_0^1 S(s) ds$$
$$+ \frac{W(1)}{2f(0)} \left[ \int_0^1 W(s) ds - \frac{W(1)}{2} \right]$$
$$=: Y.$$

Similarly, according to (2.7) in Theorem 2.1, we have

$$\lim_{\beta \downarrow 0} \frac{\partial}{\partial \alpha} U(\beta, \alpha) = \lim_{\beta \downarrow 0} \left\{ \int_0^1 e^{-\beta(1-s)} dW(s) + f(0) 2\alpha \int_0^1 e^{-2\beta(1-s)} ds \right\}$$
$$= \int_0^1 dW(s) + 2\alpha f(0) \int_0^1 ds$$
$$= W(1) + 2\alpha f(0),$$

and therefore

$$\hat{\alpha}(0+) = -\frac{W(1)}{2f(0)},$$

which is same as $\hat{\alpha}(0-)$. The derivative of $U(\beta, \alpha)$ with respect to $\beta$ at zero from righthand side satisfies

$$\frac{\partial}{\partial \beta} U(\beta, \alpha) = -\int_0^1 \int_{s+}^1 e^{-\beta(t-s)} dS(t) dW(s) - \beta \int_0^1 \int_s^1 e^{-\beta(t-s)}(s-t) dS(t) dW(s)$$
$$+ \alpha \int_0^1 e^{-\beta(1-s)}(s-1) dW(s)$$
$$+ f(0) \left\{ 2\beta \int_0^1 \left( \int_s^1 e^{-\beta(t-s)} dS(t) \right)^2 ds \right.$$
$$+ \beta^2 \int_0^1 2 \left( \int_s^1 e^{-\beta(t-s)} dS(t) \right)$$
$$\times \left( \int_s^1 e^{-\beta(t-s)}(s-t) dS(t) \right) ds$$
$$+ \alpha^2 \int_0^1 e^{-2\beta(1-s)} 2(s-1) ds$$
$$- 2\alpha \int_0^1 e^{-\beta(1-s)} \left( \int_s^1 e^{-\beta(t-s)} dS(t) \right) ds$$
$$\left. - 2\alpha\beta \int_0^1 \int_s^1 e^{-\beta(1+t-2s)}(2s-t-1) dS(t) ds \right\}.$$

Taking the limit $\beta \downarrow 0$ and using the remark in Section 2, we have

$$\lim_{\beta \downarrow 0} \frac{\partial}{\partial \beta} U(\beta, \hat{\alpha}(\beta))$$
$$\to -\int_0^1 \int_{s+}^1 dS(t) dW(s) + \hat{\alpha}(0+) \int_0^1 (s-1) dW(s)$$
$$+ f(0) \left\{ \hat{\alpha}^2(0+) \int_0^1 2(s-1) ds - 2\hat{\alpha}(0+) \int_0^1 \int_s^1 dS(t) ds \right\}$$
$$= -S(1)W(1) + \int_0^1 S(s) dW(s) + 1 + \hat{\alpha}(0+) \left[ [(s-1)W(s)]_0^1 - \int_0^1 W(s) ds \right]$$
$$+ f(0) \left\{ -\hat{\alpha}^2(0+) - 2\hat{\alpha}(0+) \left[ S(1) - \int_0^1 S(s) ds \right] \right\}$$
$$= \int_0^1 S(s) dW(s) - W(1) \int_0^1 S(s) ds + \frac{W(1)}{2f(0)} \left[ \int_0^1 W(s) ds - \frac{W(1)}{2} \right] + 1$$
$$= Y + 1.$$

Therefore, the pile-up probability in (3.3) can be expressed in terms of $Y$ as

$$\lim_{n\to\infty} P(\hat{\theta}_n^{(J)} = 1) = P[Y < 0 \text{ and } Y + 1 > 0]$$
$$= P[-1 < Y < 0].$$

## 3.2. Exact likelihood estimation

In this section, we consider pile-up probabilities associated with the estimator that maximizes the exact Laplace likelihood. For $\theta \leq 1$, the joint density of $(\boldsymbol{x}_n, z_{init})$ satisfies

$$f(\boldsymbol{x}_n, z_{init}) = \prod_{t=0}^{n} f(z_t) = \left(\frac{1}{2\sigma}\right)^{n+1} \exp\left(-\frac{\sum_{t=0}^{n} |z_t|}{\sigma}\right)$$
$$= \left(\frac{1}{2\sigma}\right)^{n+1} \exp\left\{-\frac{[\ell_n(\theta, z_{init}) - \ell_n(1, Z_0)] + \ell_n(1, Z_0)}{\sigma}\right\}$$
$$= \left(\frac{1}{2\sigma}\right)^{n+1} \exp\left(-\frac{\sum_{t=0}^{n} |Z_t|}{\sigma}\right) e^{-U_n(\beta,\alpha)}.$$

Integrating out the augmented variable $z_{init}$, we obtain

$$\int_{-\infty}^{\infty} f(\boldsymbol{x}_n, z_{init}) dz_{init} = \left(\frac{1}{2\sigma}\right)^{n+1} \exp\left(-\frac{\sum_{t=0}^{n} |Z_t|}{\sigma}\right) \frac{\sigma}{\sqrt{n}} \int_{-\infty}^{\infty} e^{-U_n(\beta,\alpha)} d\alpha,$$

since under the parameterization (2.4), $dz_{init} = (\sigma/\sqrt{n})d\alpha$. The Laplace log-likelihood of $(\theta, \sigma)$ given $\boldsymbol{x}_n$ then satisfies

$$\ell_n^*(\theta, \sigma) \equiv \log \int_{-\infty}^{\infty} f(\boldsymbol{x}_n, z_{init}) dz_{init}$$
$$= -(n+1)\log(2\sigma) - \frac{\sum_{t=0}^{n} |Z_t|}{\sigma} + \log\left(\frac{\sigma}{\sqrt{n}}\right) + \log \int_{-\infty}^{\infty} e^{-U_n(\beta,\alpha)} d\alpha,$$

where the last term does not depend on $\sigma$ as $n \to \infty$. So maximizing $\ell_n^*$ with respect to $\theta \leq 1$ is approximately the same as maximizing

(3.5) $$U_n^*(\beta) = \log \int_{-\infty}^{\infty} e^{-U_n(\beta,\alpha)} d\alpha$$

with respect to $\beta \leq 0$,

Similarly, for $\theta > 1$, the Laplace log-likelihood of $(\theta, \sigma)$ is

$$\ell_n^*(\theta, \sigma) \equiv \log \int_{-\infty}^{\infty} f(\boldsymbol{x}_n, z_{init}) dz_{init}$$
$$= -n\log|\theta| - (n+1)\log(2\sigma) - \frac{\sum_{t=0}^{n} |Z_t|}{\sigma|\theta|}$$
$$+ \log\left(\frac{\sigma}{\sqrt{n}}\right) + \log \int_{-\infty}^{\infty} e^{-U_n(\beta,\alpha)|\theta|^{-1}} d\alpha,$$

where again the last term does not depend on $\sigma$ as $n \to \infty$. As above, maximizing $\ell_n^*$ with respect to $\theta > 1$ is equivalent to maximizing

(3.6) $$U_n^*(\beta) = \log \int_{-\infty}^{\infty} e^{-U_n(\beta,\alpha)n/(n+\beta)} d\alpha$$

for $\beta > 0$.

A heuristic argument based on the process convergence of $U_n$ to $U$ suggests that

$$\text{(3.7)} \qquad U_n^*(\beta) \to U^*(\beta) = \log \int_{-\infty}^{\infty} e^{-U(\beta,\alpha)}\, d\alpha,$$

where $U_n^*$ is specified by (3.5) for $\beta \leq 0$ and by (3.6) for $\beta > 0$. Now if $\hat{\beta}_n^{(E)}$ denotes the local maximum of the exact likelihood, or alternatively the maximizer of $U_n^*(\beta)$ that is closest to 0, then the convergence in (3.7) suggests convergence in distribution for the local maximizer of the exact likelihood, i.e.,

$$\text{(3.8)} \qquad n(\hat{\theta}_n^{(E)} - 1) = \hat{\beta}_n^{(E)} \xrightarrow{d} \hat{\beta}^{(E)},$$

where $\hat{\beta}^{(E)}$ is the local maximizer of $U^*(\beta)$ that is closest to 0.

The limiting pile-up probabilities for $\hat{\theta}_n^{(E)}$ are calculated from

$$\lim_{n\to\infty} P(\hat{\theta}_n^{(E)} = 1) = \lim_{n\to\infty} P(\hat{\beta}_n^{(E)} = 0) = P(\hat{\beta}^{(E)} = 0)$$
$$= P\left(\lim_{\beta\uparrow 0} \frac{\partial}{\partial \beta} U^*(\beta) > 0 \text{ and } \lim_{\beta\downarrow 0} \frac{\partial}{\partial \beta} U^*(\beta) < 0\right).$$

Fortunately, the right- and left-hand derivatives of $U^*$ can be computed explicitly. These are found to be

$$\lim_{\beta\uparrow 0} \frac{\partial}{\partial \beta} U^*(\beta) = -\frac{W^2(1)}{4f(0)} + \frac{W(1)}{2f(0)} \int_0^1 W(s)ds - W(1) \int_0^1 S(s)ds + \int_0^1 S(s)dW(s)$$
$$+ \frac{1}{2}$$
$$= Y + \frac{1}{2},$$
$$\lim_{\beta\downarrow 0} \frac{\partial}{\partial \beta} U^*(\beta) = -\frac{W^2(1)}{4f(0)} + \frac{W(1)}{2f(0)} \int_0^1 W(s)ds - W(1) \int_0^1 S(s)ds + \int_0^1 S(s)dW(s)$$
$$+ \frac{1}{2}$$
$$= Y + \frac{1}{2},$$

where $Y$ is defined in (3.4). The limiting pile-up probability for $\hat{\theta}_n^{(E)}$ is then

$$\lim_{n\to\infty} P(\hat{\theta}_n^{(E)} = 1) = P\left[-\frac{1}{2} < Y < -\frac{1}{2}\right] = 0.$$

### 3.3. Remarks

Here we collect several remarks concerning the results of Sections 3.1 and 3.2.

**Remark 1.** Under the assumptions of Theorem 2.1, the asymptotic pile-up probability for estimator $\hat{\theta}_n^{(J)}$ based on the joint likelihood is always positive. On the other hand, the asymptotic pile-up probability for estimator $\hat{\theta}_n^{(E)}$ based on the exact likelihood is zero.

**Remark 2.** The two estimators of $\theta_0$ considered in Sections 3.1 and 3.2 were defined as the local optimizers of objective functions that were closest to 1. One could also consider the global optimizers of these objective functions. For example, the exact MLE in the Gaussian case was considered in Davis and Dunsmuir [6] and Davis, Chen and Dunsmuir [5] and has a different limiting distribution than the local MLE. In our case, there will be a positive asymptotic pile-up probability for the global maximum of the joint likelihood and a zero asymptotic pile-up probability for the global maximum of the exact likelihood.

**Remark 3.** Suppose $Z_t$ has a Laplace distribution with the density function

$$f_Z(z) = \frac{1}{2\sigma} e^{-|z|/\sigma}.$$

Then $Y$ defined in (3.4) satisfies

$$(3.9) \qquad Y = \int_0^1 [W(1)s - W(s)]\, dV(s) - \frac{1}{2},$$

where $W(s)$ and $V(s)$ are independent standard Brownian motions. To prove (3.9), note that the constant $c$ in (2.14) is equal to 1 so that

$$S(t) = W(t) + V(t).$$

In the following calculations, we use the well-known Itô formula

$$\int_0^1 W(s) dW(s) = \frac{W^2(1)}{2} - \frac{1}{2}.$$

Since $f(0) = 1/2$, the random variable $Y$ defined in (3.4) can be further simplified in terms of $W(t)$ and $V(t)$ as

$$\begin{aligned}
Y &= \int_0^1 S(s)dW(s) - W(1)\int_0^1 S(s)ds + \frac{W(1)}{2f(0)}\left[\int_0^1 W(s)ds - \frac{W(1)}{2}\right]\\
&= \int_0^1 V(s)dW(s) + \int_0^1 W(s)dW(s) - W(1)\int_0^1 V(s)ds - W(1)\int_0^1 W(s)ds \\
&\quad + W(1)\int_0^1 W(s)ds - \frac{W^2(1)}{2}\\
&= V(1)W(1) - \int_0^1 W(s)dV(s) + \frac{W^2(1)}{2} - \frac{1}{2} - W(1)\left[V(1) - \int_0^1 s\,dV(s)\right]\\
&\quad - \frac{W^2(1)}{2}\\
&= \int_0^1 [W(1)s - W(s)]\, dV(s) - \frac{1}{2}.
\end{aligned}$$

Therefore, the pile-up probability for Laplace innovations is

$$\begin{aligned}
&P(-1 < Y < 0) \\
&= P\left(-\frac{1}{2} < \int_0^1 [W(1)s - W(s)]\,dV(s) < \frac{1}{2}\right) \\
&= E\left[P\left(-\frac{1}{2} < \int_0^1 [W(1)s - W(s)]\,dV(s) < \frac{1}{2}\right) \bigg| W(t) \text{ on } t \in [0,1]\right] \\
&= E\left[P\left(-\frac{1}{2}\left\{\int_0^1 [W(1)s - W(s)]^2 ds\right\}^{-1/2} < U \right.\right. \\
&\qquad\qquad \left.\left. < \frac{1}{2}\left\{\int_0^1 [W(1)s - W(s)]^2 ds\right\}^{-1/2}\right)\right] \\
&= E\left[\Phi\left(\frac{1}{2}\left\{\int_0^1 [W(1)s - W(s)]^2 ds\right\}^{-1/2}\right)\right. \\
&\qquad\qquad \left. - \Phi\left(-\frac{1}{2}\left\{\int_0^1 [W(1)s - W(s)]^2 ds\right\}^{-1/2}\right)\right] \\
&\approx 0.820,
\end{aligned}$$

where $U$ has the standard normal distribution and $\Phi(\cdot)$ is the corresponding cumulative distribution function. This pile-up probability, which was computed via simulation based on 100000 replications of $W(t)$ on $[0, 1]$, has a standard error of 0.0010.

**Remark 4.** From the limiting result (3.2), it follows that the random variable $Z_0$ can be *estimated* consistently. It may seem odd to have a consistent estimate of a noise term in a moving average process. On the other hand, an MA(1) process with a unit root is both invertible and non-invertible. That is, $Z_0$ is an element of the two Hilbert spaces generated by the linear span of $\{X_t, t \leq 0\}$ and $\{X_t, t \geq 1\}$, respectively. It is the latter Hilbert space which allows for consistent estimation of $Z_0$.

## 4. Numerical simulation

In this section, we compute the asymptotic pile-up probabilities associated with the estimator $\hat{\theta}^{(J)}$ which maximizes the joint Laplace likelihood for several different noise distributions. The empirical properties of estimators $\hat{\theta}_n^{(J)}$ (the local maximizer of the joint Laplace likelihood) and $\hat{\theta}_n^{(E)}$ (the local maximizer of the exact Laplace likelihood) for finite samples are compared with each other and with the corresponding asymptotic theory.

For approximating the asymptotic pile-up probabilities and limiting distribution of $\hat{\beta}_n^{(J)}$, we first simulate 100000 replications of independent standard Wiener processes $W(t)$ and $V(t)$ on $[0, 1]$ in which $W(t)$ and $V(t)$ are approximated by the partial sums $W(t) = \sum_{j=1}^{[10000t]} W_j/\sqrt{10000}$ and $V(t) = \sum_{j=1}^{[10000t]} V_j/\sqrt{10000}$, where $\{W_j\}$ and $\{V_j\}$ are independent standard normal random variables. From the simulation of $W(t)$ and $V(t)$, the distribution of the limit random variable $\hat{\beta}^{(J)}$ can be tabulated and the pile-up probability $P(-1 < Y < 0)$ estimated, where $Y$ is given in (3.4). The empirical pile-up probabilities and their asymptotic limits are displayed in Table 1 for different noise distributions: Laplace, Gaussian, uniform, and $t$ with 5 degrees of freedom. Notice that there is good agreement between the asymptotic and empirical probabilities for sample sizes as small as 50.

For examining the empirical performance of the local maximizers $\hat{\theta}_n^{(J)}$ and $\hat{\theta}_n^{(E)}$, we only consider the process generated with Laplace noise with $\sigma = 1$ and sample sizes $n = 20, 50, 100, 200$. For each setup, 1000 realizations of the MA(1) process with $\theta_0 = 1$ are generated and the estimates $\hat{\theta}_n^{(J)}$ and $\hat{\theta}_n^{(E)}$ and their corresponding estimates of the scale parameter are obtained. The estimation results are summarized in Table 2. For comparison, the standard deviation based on the limit distributions of $\hat{\theta}_n^{(J)}$ and $\hat{\theta}_n^{(E)}$ are also reported (denoted by asymp in the table), which are obtained numerically based on 100000 replicates of the limit process $U$. Generally speaking, the empirical root mean square errors are very close to their asymptotic values even for very small samples. Moreover, the estimation error of $\hat{\theta}_n^{(J)}$ is about 1/2 the estimation error of $\hat{\theta}_n^{(E)}$, which indicates the superiority of using the joint likelihood over exact likelihood when $\theta_0 = 1$.

We also considered performance of the two estimators $\hat{\theta}_n^{(J)}$ and $\hat{\theta}_n^{(E)}$ in the case when $\theta_0 \neq 1$. A limit theory for these estimators can be derived in this case by assuming that the true value $\theta_0$ is near 1. That is, we can parameterize the MA(1) parameter by $\theta_0 = 1 + \gamma/n$ (e.g., Davis and Dunsmuir [6]). While we have not pursued the theory in the near unit root case, the relative performance of these estimators was compared in a limited simulation study. We considered 3 values of $\theta_0 = 0.8, 0.9, 0.95$ and their reciprocals $1/0.8, 1/0.9, 1/0.95$. The latter 3 cases correspond to purely non-invertible models. The results reported in Table 3 are based on the global optimization of the joint and exact likelihoods. The first two columns contain the number of realizations out of 1000 in which the estimator was invertible ($< 1$) and on the unit circle ($= 1$), respectively. For example, in the $\theta_0 = 0.8$ and $\hat{\theta}_n^{(J)}$ case, 78.9% of the realizations produced invertible models, and the empirical pile-up probability is 0.095. On the other hand, for $\theta_0 = 1/0.8$, 79.5% of the realizations produced a purely non-invertible model with an empirical pile-up probability of 0.109. Both objective functions do a reasonably good job of discriminating between invertible and non-invertible models, with a performance edge going to the exact likelihood. In terms of root mean square error, the performance of $\hat{\theta}_n^{(E)}$ is superior to $\hat{\theta}_n^{(J)}$ as $\theta_0$ moves away from the unit circle.

TABLE 1
*Empirical pile-up probabilities of the local maximizer $\hat{\theta}_n^{(J)}$ of the joint Laplace likelihood for an MA(1) with $\theta_0 = 1$ and sample sizes $n = 20, 50, 100, 200$ (based on 1000 replicates) and their asymptotic values under various noise distributions.*

| $n$ | Gau | Lap | Unif | $t(5)$ |
|---|---|---|---|---|
| 20 | 0.827 | 0.796 | 0.831 | 0.796 |
| 50 | 0.859 | 0.806 | 0.864 | 0.823 |
| 100 | 0.873 | 0.819 | 0.864 | 0.817 |
| 200 | 0.844 | 0.819 | 0.843 | 0.831 |
| 500 | 0.855 | 0.809 | 0.841 | 0.846 |
| $\infty$ | 0.873 | 0.820 | 0.862 | 0.836 |

TABLE 2
*Bias, standard deviation and root mean square error of the local maximizers $\hat{\theta}_n^{(J)}$ and $\hat{\theta}_n^{(E)}$ of the joint and exact Laplace likelihoods, respectively, for an MA(1) process generated by Laplace noise with $\theta_0 = 1$ and $\sigma = 1$ (1000 replications).*

| $n$ | | $\hat{\theta}_n^{(J)}$ | $\hat{\theta}_n^{(E)}$ |
|---|---|---|---|
| $n = 20$ | bias | -0.003 | -0.006 |
| | s.d. | 0.066 | 0.144 |
| | rmse | 0.066 | 0.144 |
| | asymp | 0.053 | 0.121 |
| $n = 50$ | bias | -0.000 | 0.000 |
| | s.d. | 0.021 | 0.057 |
| | rmse | 0.021 | 0.057 |
| | asymp | 0.021 | 0.048 |
| $n = 100$ | bias | -0.000 | 0.001 |
| | s.d. | 0.011 | 0.030 |
| | rmse | 0.011 | 0.030 |
| | asymp | 0.011 | 0.024 |
| $n = 200$ | bias | 0.000 | 0.001 |
| | s.d. | 0.006 | 0.014 |
| | rmse | 0.006 | 0.014 |
| | asymp | 0.005 | 0.012 |

TABLE 3
*Bias, standard deviation and root mean square error of the global maximizers $\hat{\theta}_n^{(J)}$ and $\hat{\theta}_n^{(E)}$ of the joint and exact Laplace likelihoods, respectively, for an MA(1) process generated by Laplace noise with $\theta_0 = 0.8, 0.9, 0.95, 1/0.95, 1/0.9, 1/0.8$, $\sigma = 1$, and $n = 50$ based on 1000 replications. First 2 columns record the number of times (out of 1000) that the estimates were less than 1 (invertible) and equal to 1 (unit root).*

| $\theta_0$ |  | $< 1$ | $= 1$ | bias | s.d. | rmse |
|---|---|---|---|---|---|---|
| 0.8 | $\hat{\theta}_{50}^{(J)}$ | 789 | 95 | 0.0734 | 0.1973 | 0.2105 |
|  | $\hat{\theta}_{50}^{(E)}$ | 873 | 19 | 0.0498 | 0.1753 | 0.1822 |
| 0.9 | $\hat{\theta}_{50}^{(J)}$ | 557 | 322 | 0.0578 | 0.1398 | 0.1513 |
|  | $\hat{\theta}_{50}^{(E)}$ | 767 | 93 | 0.0327 | 0.0933 | 0.0989 |
| 0.95 | $\hat{\theta}_{50}^{(J)}$ | 404 | 503 | 0.0322 | 0.0708 | 0.0778 |
|  | $\hat{\theta}_{50}^{(E)}$ | 632 | 168 | 0.0235 | 0.0821 | 0.0854 |
| 1/0.95 | $\hat{\theta}_{50}^{(J)}$ | 90 | 540 | -0.0315 | 0.0763 | 0.0825 |
|  | $\hat{\theta}_{50}^{(E)}$ | 286 | 114 | -0.0207 | 0.0890 | 0.0914 |
| 1/0.9 | $\hat{\theta}_{50}^{(J)}$ | 89 | 299 | -0.0389 | 0.1227 | 0.1287 |
|  | $\hat{\theta}_{50}^{(E)}$ | 207 | 71 | -0.0327 | 0.1218 | 0.1261 |
| 1/0.8 | $\hat{\theta}_{50}^{(J)}$ | 96 | 109 | -0.0338 | 0.2645 | 0.2666 |
|  | $\hat{\theta}_{50}^{(E)}$ | 149 | 19 | -0.0492 | 0.2280 | 0.2333 |

**Remark.** The LAD estimate of $\theta_0$ is obtained by minimizing the objective function given in (2.2) with $z_{init} = 0$. Although we have not considered the asymptotic pile-up in this case, the estimator does not perform as well as $\hat{\theta}_n^{(J)}$ and $\hat{\theta}_n^{(E)}$. For example, in simulation results, not reported here, the rmse of the LAD estimator tended to be twice as large as the rmse for the exact MLE.